\documentclass[11pt]{amsart}
\usepackage[utf8]{inputenc}
\usepackage[margin=1.2in]{geometry}

\usepackage{amsmath}
\usepackage{amsfonts}
\usepackage{amssymb}
\usepackage{amsthm}
\usepackage{mathtools}
\usepackage{caption}
\usepackage{subcaption}
\usepackage{bbm}
\usepackage[export]{adjustbox}
\usepackage{stmaryrd}
\usepackage[all]{xy}

\usepackage{tikz-cd}
\usetikzlibrary{matrix,arrows.meta,calc}
\usepackage{graphicx}
\usepackage{float}
\usepackage{epstopdf}
\usepackage{overpic}
\usepackage{microtype}
\usepackage{placeins}

\usepackage[linktocpage,hypertexnames=false]{hyperref}
\hypersetup{
    colorlinks=true,
    linkcolor=blue,
    citecolor=blue,
    urlcolor=blue,
    pdftitle={Nearby Lagrangians in the cotangent bundle of the Klein bottle},
    pdfsubject={Neck stretching, invariant sphere rulings, and relative reconstruction}
}

\usepackage{color}
\definecolor{note}{rgb}{0,0,1}
\definecolor{curveA}{RGB}{35,91,135}
\definecolor{curveB}{RGB}{164,76,34}
\newtheorem{theorem}{Theorem}[section]
\newtheorem{proposition}[theorem]{Proposition}
\newtheorem{lemma}[theorem]{Lemma}
\theoremstyle{definition}

\theoremstyle{remark}
\newtheorem{remark}[theorem]{Remark}
\numberwithin{equation}{section}

\DeclareMathOperator{\Ham}{Ham}

\DeclareMathOperator{\id}{id}
\newcommand{\R}{\mathbb R}
\newcommand{\Z}{\mathbb Z}
\newcommand{\C}{\mathbb C}

\newcommand{\op}{\operatorname}

\newcommand{\hs}{\vskip.1in}

\title[Nearby Lagrangians for the Klein bottle]{Nearby Lagrangians in the cotangent bundle of the Klein bottle}
\author{Yuan Gao}
\address{School of Mathematics, Nanjing University, Nanjing, Jiangsu Province 210093, China}
\email{yuangao@nju.edu.cn}
\urladdr{}
\author{Weiwei Wu}
\address{School of Mathematical Sciences, Zhejiang University, Hangzhou, Zhejiang Province 310027, China}
\email{wwwu9693@zju.edu.cn}
\urladdr{}
\author{Tianyu Yuan}
\address{School of Mathematical Sciences, Eastern Institute of Technology, Ningbo, Zhejiang, 315200, China}
\email{tyyuan@eitech.edu.cn}
\urladdr{}
\date{\today}

\subjclass[2020]{Primary 53D12; Secondary 53D35, 53D42}
\keywords{Nearby Lagrangian conjecture, Klein bottle, equivariant symplectic isotopy}

\begin{document}

\begin{abstract}
    We prove the nearby Lagrangian conjecture for the Klein bottle: every closed connected embedded exact Lagrangian submanifold of its cotangent bundle is Hamiltonian isotopic to the zero section by a compactly supported isotopy.
    The proof passes to the torus double cover and builds on the theorem of Dimitroglou Rizell, Goodman, and Ivrii for the two-dimensional torus.
    The main step is to construct an isotopy compatible with the covering involution.
    We achieve this using holomorphic sphere rulings in $S^2\times S^2$ adapted to the zero section, allowing the construction to descend to the Klein bottle.
\end{abstract}

\maketitle

\tableofcontents

\section{Introduction}\label{section-intro}

The nearby Lagrangian conjecture asks whether a closed embedded exact Lagrangian in a cotangent bundle of a closed manifold is Hamiltonian isotopic to its zero section.
In this paper, we prove the conjecture in the case of the Klein bottle.

\begin{theorem}\label{thm-main}
    Let $K$ be the Klein bottle and let $(T^*K, \omega_K=d\lambda_K)$ be its cotangent bundle.
    Every closed connected embedded exact Lagrangian submanifold of $(T^*K, d\lambda_K)$ is Hamiltonian isotopic to the zero section by a smooth isotopy supported in one compact subset of $T^*K$.
\end{theorem}

Earlier work establishes topological and categorical constraints on closed exact Lagrangians in cotangent bundles.
For a simply connected spin base, Fukaya, Seidel, and Smith \cite{FSS2008} proved that the cotangent projection of a closed exact spin Lagrangian with vanishing Maslov class induces an isomorphism on complex cohomology.
Under the same hypotheses, Nadler \cite{Nadler2009} identified the Lagrangian with the zero section in the Fukaya category over $\mathbb C$, up to a grading shift.
The work of Abouzaid \cite{Abouzaid2012Homotopy}, Kragh \cite{Kragh2013}, and Abouzaid--Kragh \cite{AK} shows that the cotangent projection of every closed exact Lagrangian is a simple homotopy equivalence.

For closed surfaces, Hind proved the conjecture for $S^2$ \cite{Hind2004}, and Li and Wu gave an alternative proof using symplectic cutting \cite[Theorem~6.6]{LiWu2012}.
In $T^*\mathbb{RP}^2$, Hind proved Hamiltonian uniqueness for Lagrangian projective planes homotopic to the zero section \cite[Corollary~8]{Hind2012}.
Hind, Pinsonnault, and Wu proved that the space of Lagrangian projective planes in $T^*\mathbb{RP}^2$ is weakly contractible \cite[Theorem~1.3]{HPW2016}, and Adaloglou gave a compactification proof of Hamiltonian uniqueness \cite{Adaloglou2025}.
Dimitroglou Rizell, Goodman, and Ivrii established the conjecture for $T^2$ \cite[Theorem~B]{DRGI}, providing the starting point of our argument.

For the Klein bottle, Atallah, Chass\'e, Leclercq, and Shelukhin proved that every Lagrangian Klein bottle in $T^*K$ is $H$-exact, meaning that symplectic area vanishes on relative two-dimensional homology classes \cite[Theorem~44]{ACLS2025}.
This condition is weaker than exactness of the restricted Liouville form, which is assumed in Theorem~\ref{thm-main}.
Their argument also uses the torus double cover.

\subsection*{Outline of the proof}

The main idea of this paper is to pass the problem of the Klein bottle to its double cover.
Let $L\subset T^*K$ be as in Theorem~\ref{thm-main}.
The Abouzaid--Kragh theorem implies that its full lift $E\subset T^*T^2$ is a connected exact torus.
Write $\sigma: T^*T^2 \to T^*K$ for the covering involution and $Z$ for the zero section of $T^*T^2$.
We first use the torus theorem of Dimitroglou Rizell, Goodman, and Ivrii \cite[Theorem~B]{DRGI} to obtain a compactly supported Hamiltonian diffeomorphism $\psi$ taking $E\subset T^*T^2$ to the zero section $Z$ of $T^*T^2$.
Descent to $T^*K$ requires an isotopy commuting with $\sigma$.
As a first step, Lemma~\ref{lemma-normalization} allows us to choose $\psi$ so that the conjugated involution $\rho=\psi\sigma\psi^{-1}$ agrees with $\sigma$ near $Z$.

To construct this path, we consider the equal-area compactification $Y=S^2\times S^2$ of a square cotangent domain in $T^*T^2$, where $Z$ becomes the standard Clifford torus.
Gromov's theorem gives two transverse sphere rulings for each compatible almost complex structure on $Y$ \cite[Section~2.4.A1]{Gromov}.
Our geometric input is an equivariant adaptation of the neck-stretching and smoothing constructions in \cite[Sections~4--5]{DRGI}.
We construct a single $\rho$-invariant compatible almost complex structure whose two rulings contain circle families of spheres meeting $Z$ in its two families of coordinate circles (Proposition~\ref{prop-marking}).
Intersections of ruling leaves \cite[pp.~11--13]{Evans}, followed by the equivariant annular straightening of Lemma~\ref{lemma-annulus}, then give a smooth path of diffeomorphisms $f_t:Y\to Y$ satisfying
\[
    f_0=\psi^{-1},\qquad f_t\rho=\sigma f_t,\qquad f_1(Z)=Z.
\]

Using the convex-form argument of \cite[pp.~12--13]{Evans}, we apply a relative Moser deformation to make this path symplectic while preserving its initial value, conjugacy relation, and zero-section endpoint condition.
A Hamiltonian correction near the caps then gives an equivariant symplectic isotopy with common compact support, whose descent to $T^*K$ is Hamiltonian because $H_c^1(T^*K;\R)=0$.

\subsection*{Structure of the paper}

In Section~\ref{section-cover}, we normalize the involution near the zero section, state the central equivariant symplectic isotopy result, and introduce the square compactification together with the extended involution and intersection properties.
In Section~\ref{section-planes}, we construct the two rulings with the prescribed intersections with the zero section.
In Section~\ref{section-reconstruction}, we use the two rulings to get the desired smooth isotopy and apply Moser deformation to make it symplectic.
Finally, we prove the main theorem in Section~\ref{section-descent}.

\subsection*{Acknowledgments}

The authors thank Nik Adaloglou, Johannes Hauber, and Jialiang Chen for helpful discussions.
The authors formulated the research problem, proposed the double-cover approach, and outlined a program using an equivariant adaptation of the foliation method in \cite{DRGI}.
They supplied concrete mathematical directions for constructing and using the foliations in the compactified setting.
Under the authors' direction, a generative large language model produced a preliminary draft of more than 40 pages.
The authors guided the mathematical development throughout, examined and revised the generated arguments, and extensively revised the manuscript for clarity and precision.
They take full responsibility for the mathematical content of the final manuscript.

\section{The orientation cover and compactification}\label{section-cover}

\subsection{The orientation cover and normalization}

Let $Z$ denote the zero section of $T^*T^2$.
We use coordinates
\begin{equation*}
    (x,y,p_x,p_y)\in (\R/\Z)^2\times\R^2,\qquad \lambda=p_x dx+p_y dy,\qquad \omega=d\lambda.
\end{equation*}
Write $q=(x,y)$ and $p=(p_x,p_y)$ for the base and fiber coordinates.
The cotangent lift of the orientation-cover involution is
\begin{equation}\label{eq-sigma}
    \sigma(x,y,p_x,p_y)=(x+\frac{1}{2},-y,p_x,-p_y).
\end{equation}
Therefore, $T^*T^2/\langle\sigma\rangle=T^*K$.
Throughout the paper, an isotopy is called \emph{compactly supported} if all its maps are the identity outside one common fixed compact set.
We use the convention $\iota_{X_H}\omega=-dH$.
We write $\Ham_c(T^*T^2)$ for the group of time-one maps of compactly supported Hamiltonian isotopies.

The key technical part of the paper is to construct the following equivariant symplectic isotopy:
\begin{proposition}\label{prop-equivariant}
    Let $\psi\in\Ham_c(T^*T^2)$, and suppose that $\rho=\psi\sigma\psi^{-1}$ agrees with $\sigma$ on a neighborhood of $Z$.
    There is a smooth symplectic isotopy $h_t$, supported in one common compact subset of $T^*T^2$, such that
    \begin{equation*}
        h_0=\id,\qquad h_t\sigma=\sigma h_t\quad (0\leq t\leq1),
        \qquad h_1 (\psi^{-1}(Z))=Z.
    \end{equation*}
\end{proposition}

We will prove the proposition in Sections~\ref{section-planes}--\ref{section-reconstruction}.
We first show that any closed exact Lagrangian of $T^*K$ gives rise to such a $\psi$ that satisfies the assumptions of Proposition \ref{prop-equivariant}.

\begin{lemma}\label{lemma-normalization}
    If $L\subset T^*K$ is a closed connected embedded exact Lagrangian, then its inverse image $E\subset T^*T^2$ is a connected exact torus.
    Moreover, there is a $\psi\in\Ham_c(T^*T^2)$ with
    \begin{equation}\label{eq-normalization}
        \psi(E)=Z,\qquad \psi\sigma\psi^{-1}=\sigma\quad \text{near } Z.
    \end{equation}
\end{lemma}
\begin{proof}
    \emph{Step 1: Lifting the Lagrangian.}
    The cotangent projection $L\to K$ is a homotopy equivalence by \cite[Theorem~1]{AK}.
    Its induced map on fundamental groups therefore pulls the orientation character of $K$ back to a surjection onto $\Z/2$.
    The associated cover $E\to L$ is connected and has fundamental group $\Z^2$.
    The classification of closed surfaces then identifies $E$ with a torus.
    Pulling back a primitive for $\lambda_K|_L$ proves exactness.
    By \cite[Theorem~B]{DRGI}, there is a $\phi\in\Ham_c(T^*T^2)$ with $\phi(E)=Z$.

    \hs
    \noindent\emph{Step 2: Normalizing the action on the zero section.}
    Let $\rho_0=\phi\sigma\phi^{-1}$.
    This involution acts freely on $Z$, and
    \begin{equation*}
        (\rho_0|_Z)_*=
        \begin{pmatrix}
            1 & 0 \\
            0 & -1
        \end{pmatrix}
        \quad \text{on } H_1(Z;\Z),
    \end{equation*}
    since $Z\hookrightarrow T^*T^2$ is a homotopy equivalence and $\phi$ is isotopic to the identity.
    The quotient of $Z$ by $\rho_0$ is a Klein bottle, and its orientation cover gives a diffeomorphism $c:Z\to Z$ which conjugates $\rho_0|_Z$ to $\sigma|_Z$.
    The matrix $c_*$ commutes with $\operatorname{diag}(1,-1)$, so it is diagonal with entries $\pm1$.
    The corresponding linear torus map $d$ commutes with $\sigma|_Z$, including its half-translation.
    Replacing $c$ by $d^{-1}c$ then makes $c_*$ the identity.
    It follows from the mapping class group calculation for the torus \cite[Theorem~2.5]{FM} that $c$ is isotopic to the identity.

    Let $c_t$ be an isotopy from $\id_Z$ to $c$, and let $v_t=\dot c_t\circ c_t^{-1}$ be its generating vector field.
    The cotangent lift of $c_t$ is generated by $H_t(q,p)=p(v_t(q))$.
    Cutting off $H_t$ in the fibers gives a compactly supported Hamiltonian isotopy $\kappa_t$ preserving $Z$ and satisfying $\kappa_t|_Z=c_t$.
    Define
    \[
        \rho_1=\kappa_1\rho_0\kappa_1^{-1}.
    \]
    Since $\kappa_1|_Z=c$, we have $\rho_1|_Z=\sigma|_Z$.

    \hs
    \noindent\emph{Step 3: Normalizing the germ along the zero section.}
    By the equivariant coisotropic neighborhood theorem \cite[Section~30.2, footnote~19]{CdS}, the identity map of $Z$ extends to a germ of symplectomorphisms $g$ of a neighborhood of $Z$, satisfying
    \[
        g|_Z=\id_Z, \qquad
        g\rho_1g^{-1}=\sigma
        \quad\text{near }Z.
    \]

    We realize $g$ by a compactly supported Hamiltonian diffeomorphism.
    Let $\delta_s(q,p)=(q,sp)$.
    Since $g$ is symplectic and fixes $Z$ pointwise,
    \[
        Dg|_{(q,0)} =
        \begin{pmatrix}
            I & B(q) \\
            0 & I
        \end{pmatrix},
        \qquad B(q)=B(q)^T.
    \]
    Taylor's formula therefore shows that the germs of symplectomorphisms
    \[
        g_s=\delta_s^{-1}\circ g\circ\delta_s,
        \qquad 0<s\leq1,
    \]
    extend smoothly to $g_0=\id$.
    On a common neighborhood of $Z$, set
    \[
        V_s=\partial_sg_s\circ g_s^{-1}.
    \]
    Since $g_s$ is symplectic and fixes $Z$ pointwise, the one-form $\iota_{V_s}\omega$ is closed and vanishes along $Z$.
    This implies $V_s$ is generated by a family of Hamiltonians $-F_s$ satisfying $F_s|_Z=0$.

    Cutting off the Hamiltonians $-F_s$ by a fixed cutoff function equal to one near $Z$, we obtain $g_s$ on a smaller neighborhood for all $s$.
    Its time-one map $G\in\Ham_c(T^*T^2)$ therefore fixes $Z$ pointwise and agrees with $g$ near $Z$.

    Finally, set $\psi=G\circ\kappa\circ\phi$.
    Then $\psi\in\Ham_c(T^*T^2)$, $\psi(E)=Z$, and
    \[
        \psi\sigma\psi^{-1}=G\rho_1 G^{-1}=\sigma
        \quad\text{near }Z.
    \]
    This proves \eqref{eq-normalization}.
\end{proof}

\subsection{The square compactification}\label{section-square}

Fix $\psi$ and $\rho$ as in Proposition~\ref{prop-equivariant}.
Choose $R>0$ so that $\op{supp}\psi$ lies in
\begin{equation*}
    U=(\R/\Z)^2\times (-R,R)^2.
\end{equation*}
Collapsing the boundary circles of each annular factor compactifies $U$ to
\begin{equation*}
    (Y,\omega)=S^2(2R)\times S^2(2R),
\end{equation*}
where each sphere has symplectic area $2R$.
Let $a_\pm$ be the poles of the first sphere factor, corresponding to $p_x=\pm R$, and let $b_\pm$ be the poles of the second, corresponding to $p_y=\pm R$.
The four \textit{caps} are
\[
    D_1^\pm=\{a_\pm\}\times S^2,\qquad
    D_2^\pm=S^2\times\{b_\pm\},
\]
and we set $D=\bigcup_{i,\pm}D_i^\pm$.
Their four intersection points $(a_\pm,b_\pm)$ are called \textit{the corners}.

The involution in \eqref{eq-sigma} extends to $Y$ as $\sigma=\sigma_1\times\sigma_2$, where in the annular coordinates
\[
    \sigma_1(x,p_x)=(x+\tfrac12,p_x),\qquad
    \sigma_2(y,p_y)=(-y,-p_y).
\]
On the compactified spheres, these are rotations by $\pi$ about the polar axis of the first factor and an equatorial axis of the second factor, respectively.
Thus $\sigma_1$ fixes the poles $a_\pm$, whereas $\sigma_2$ exchanges the poles $b_\pm$ and fixes exactly two equatorial points, given by $(y,p_y)=(0,0)$ and $(\tfrac12,0)$.
The map $\psi$ extends as the identity near $D$, so $\rho=\psi\sigma\psi^{-1}$ also extends to $Y$ and agrees with $\sigma$ near $D$.
Note that neither $\rho$ nor $\sigma$ fixes any corner: both exchange $(a_\pm,b_+)$ and $(a_\pm,b_-)$.

Let
\begin{equation*}
    A=[S^2\times\{\mathrm{pt}\}],\qquad B=[\{\mathrm{pt}\}\times S^2].
\end{equation*}
The caps $D_1^\pm$ represent $B$, and $D_2^\pm$ represent $A$, and $Z$ is the product of the equators.
Fix a $\sigma$-invariant $\omega$-compatible product almost complex structure $J_{\mathrm{std}}$ satisfying $J_{\mathrm{std}}\partial_x=-\partial_{p_x}$ and $J_{\mathrm{std}}\partial_y=-\partial_{p_y}$ near $Z$.
An almost complex structure on $Y$ is called \emph{standard near a set} if it agrees with $J_{\mathrm{std}}$ on a neighborhood of that set.

Fix a point $p\in Z$, let $H_i^\pm$ be the four standard symplectic disks with boundary on $Z$, whose classes form a standard basis of $H_2(Y, Z;\Z)$, and $H_i^\pm\cap D_i^\pm=1$.
Each of the disks has area $R$ and Maslov index $2$.
Let $e_1=[x]$ and $e_2=[y]$ be the positively oriented coordinate circle classes in $H_1(Z;\Z)$.
Then
\begin{equation}\label{eq-boundary}
    \partial H_i^+=-e_i,\qquad \partial H_i^-=e_i.
\end{equation}

\begin{figure}[htbp]
    \centering
    \includegraphics[width=5cm]{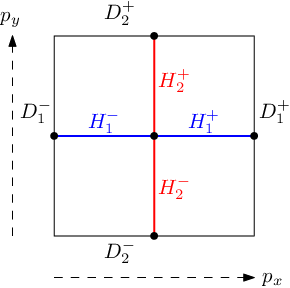}
    \caption{Momentum projection of the square compactification.
        Edge preimages are the indicated cap spheres, and the center has fiber $Z\cong T^2$.
        Each segment is the projection of $H_i^\pm$, with the complementary base coordinate fixed.
        The indicated circle collapses along an edge and both circles collapse at a corner.}
    \label{fig-cap-square}
\end{figure}

Figure~\ref{fig-cap-square} explains the cap and disk labels.

For every $\omega$-compatible almost complex structure $J$ on $Y$, the $J$-holomorphic spheres in $A$ and $B$ form two smooth transverse sphere rulings.
Each point of $Y$ lies on exactly one leaf of each ruling, and leaves in different classes intersect once transversely.
These rulings also vary smoothly with a smooth family of $J$.
If $J$ is standard near $D$, the four caps are ruling leaves.

\section{Neck stretching argument}\label{section-planes}

The goal of this section is to construct a $\rho$-invariant almost complex structure for which the two sphere rulings have prescribed intersections with $Z$.
This will be used in the reconstruction of equivariant symplectic isotopy in Section~\ref{section-reconstruction}.

\begin{proposition}\label{prop-marking}
    There is a smooth $\rho$-invariant $\omega$-compatible almost complex structure $J_1$ on $Y$, standard near $D\cup Z$, with the following property.
    For every $c,d\in\R/\Z$, the circles
    \begin{equation*}
        \{(x,c,0,0): x\in\R/\Z\},\qquad \{(d,y,0,0): y\in\R/\Z\}
    \end{equation*}
    are, respectively, the intersections of $Z$ with a $J_1$-holomorphic sphere in class $A$ and a $J_1$-holomorphic sphere in class $B$.
\end{proposition}

We stretch a neck along $Z$ and use the area calculation in Subsection~\ref{section-square} to determine the limiting configurations of ruling spheres through $Z$.
Their exterior planes form four circle families, which we straighten equivariantly and join in opposite pairs across $Z$ to symplectic spheres, following \cite[Sections~5.3--5.4]{DRGI} We then choose a compatible almost complex structure making both sphere families holomorphic.
Their coordinate-circle intersections identify $Z$ with a product of circles in the two leaf spaces, which will make the endpoint of the reconstruction in Section~\ref{subsection-leaf-maps} preserve $Z$.

\subsection{Limiting spheres and the four caps}\label{subsection-splitting}

Equip $Z=T^2$ with the flat metric $dx^2+dy^2$, and write $r=|p|$ in the cotangent coordinates.
Choose a $\rho$-invariant Weinstein neighborhood $U$ of $Z$, disjoint from $D$, on which $\rho=\sigma$, and a cosphere bundle $\Sigma\subset U$.
We stretch a neck along $\Sigma$ using the construction of \cite[Section~2.5]{DRGI}.
Attaching cylindrical ends to the two sides of $\Sigma$ identifies them with $T^*T^2$ and $Y\setminus Z$, respectively, with the fixed normalization of the Weinstein coordinates used in that construction.
We first specify the almost complex structures on these two spaces and their cylindrical ends.

On $S^*T^2$, write $p=(\cos\theta,\sin\theta)$; the contact form and Reeb vector field are
\begin{equation*}
    \alpha_0=\lambda|_{S^*T^2}=\cos\theta\,dx+\sin\theta\,dy, \qquad
    R_{\alpha_0}=\cos\theta\,\partial_x+\sin\theta\,\partial_y.
\end{equation*}
On $\R\times S^*T^2$, identified with $T^*T^2\setminus Z$ by $|p|=e^s$, use the cylindrical almost complex structure
\begin{equation}\label{eq-cylindrical-structure}
    J_{\mathrm{cyl}}\partial_x=-|p|\partial_{p_x},\qquad
    J_{\mathrm{cyl}}\partial_y=-|p|\partial_{p_y}.
\end{equation}
On $T^*T^2$, fix the almost complex structure of \cite[Section~4]{DRGI}, denoted here by $J_{\mathrm{in}}$:
\begin{equation}\label{eq-inner-structure}
    J_{\mathrm{in}}\partial_x=-f(|p|)\partial_{p_x},\qquad
    J_{\mathrm{in}}\partial_y=-f(|p|)\partial_{p_y},
\end{equation}
where $f>0$ is smooth and nondecreasing, equals $1$ near zero, and satisfies $f(r)=r$ for $r\gg0$.
The values on the fiber directions are determined by $J^2=-\id$.
Thus $J_{\mathrm{in}}$ is standard near the zero section and agrees with $J_{\mathrm{cyl}}$ near fiber infinity.
Formula~\eqref{eq-inner-structure} holds throughout the interpolation region as well, so the coordinate planes spanned by $\partial_x,\partial_{p_x}$ and by $\partial_y,\partial_{p_y}$ remain $J_{\mathrm{in}}$-invariant.
The cylinder classification in \cite[Lemma~4.2]{DRGI} applies to this model.

On $Y\setminus Z$, choose a $\rho$-invariant $\omega$-compatible almost complex structure $J_{\mathrm{out}}$ which is standard near $D$ and equals $J_{\mathrm{cyl}}$ on the end approaching $Z$.
Such a choice exists by the standard invariant extension construction: the prescribed local structures are $\sigma$-invariant, and $\rho=\sigma$ near $D\cup Z$.

With $J_{\mathrm{out}}$ and $J_{\mathrm{in}}$ fixed, the neck-stretching construction gives a family $J_\ell^{\mathrm{pre}}$ of $\omega$-compatible almost complex structures on $Y$, with cylindrical neck length tending to infinity as $\ell\to\infty$.
The family can be chosen $\rho$-invariant because the stretching takes place where $\rho=\sigma$ and commutes with this action.
Under the Weinstein coordinates, $J_\ell^{\mathrm{pre}}=J_{\mathrm{in}}$ on a fixed neighborhood of $Z$ for every $\ell$.
The almost complex structures also equal $J_{\mathrm{std}}$ on a fixed neighborhood of $D$, so the four caps remain holomorphic.

As $\ell\to\infty$, SFT compactness \cite[Section~2.5 and Theorem~2.2]{DRGI} and \cite{BEHWZ} gives subsequential limits of $J_\ell^{\mathrm{pre}}$-holomorphic ruling spheres in either class $A$ or $B$ as connected genus-zero holomorphic buildings in the same class.
Their components lie in the two completed sides specified above, with zero or finitely many intermediate levels in $(\R\times S^*T^2,J_{\mathrm{cyl}})$.
We call the components in $(Y\setminus Z,J_{\mathrm{out}})$ \emph{exterior components}, and those in $(T^*T^2,J_{\mathrm{in}})$ \emph{bottom components}.

A connected genus-zero holomorphic building in class $A$ or $B$ in these spaces is called a \emph{split sphere}.
Matched punctures have the same Reeb orbit and multiplicity; ordinary nodes join components in the same target level.
An exterior plane is a finite-energy $J_{\mathrm{out}}$-holomorphic map from $\C$ with one negative puncture at infinity, considered up to holomorphic reparametrization.
A cylinder has domain $\C^*$; in an intermediate level, a trivial cylinder has image $\R\times\gamma$ for a Reeb orbit $\gamma$.
Exterior punctures approach $|p|=0$, while bottom punctures approach $|p|=\infty$ in $T^*T^2$.

\begin{figure}[htbp]
    \centering
    \includegraphics[width=10cm]{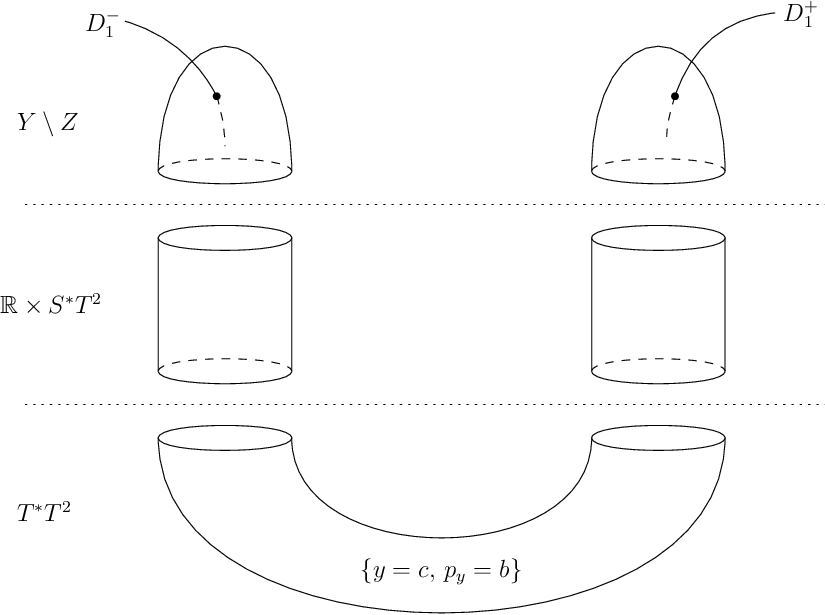}
    \caption{The class-$A$ split sphere in Lemma~\ref{lemma-splitting}.
        Each exterior plane meets the indicated cap once.
        The two columns correspond to opposite Reeb directions.
        The inner cylinder has $b=0$ in the point-constrained moduli space.
        The class-$B$ picture exchanges the coordinate directions.}
    \label{fig-building}
\end{figure}

\begin{lemma}\label{lemma-splitting}
    A split sphere in $A$ or $B$ with a nonconstant bottom component consists of two exterior planes, one primitive inner cylinder, and trivial intermediate cylinders.
    The relative classes of the two planes have area $R$.
    In class $A$, the planes meet $D_1^-$ and $D_1^+$ once, respectively, and their asymptotes are the two orientations of the same primitive $x$-geodesic.
    The bottom cylinder has image $\{y=c,\ p_y=b\}$ for some $c\in\R/\Z$ and $b\in\R$.
    For class $B$, exchange $x$ and $y$ and use the caps $D_2^\pm$.
    For a limit of embedded ruling spheres, the exterior planes are embedded.
    See Figure~\ref{fig-building}.
\end{lemma}
\begin{proof}
    We use the relative area formula to obtain the splitting without a genericity assumption on $J_{\mathrm{out}}$.
    Fix $0<\epsilon<R/3$ and a split sphere.
    Truncate its exterior ends sufficiently far down the cylindrical ends that the annuli completing their boundary curves to $Z$ have total area less than $\epsilon$.
    Each nonconstant exterior component represents a relative class (by adding a trivial cylinder with boundary on $Z$) of positive area, hence of area at least $R$ by the area formula in Subsection~\ref{section-square}.
    Its truncated area is therefore at least $R-\epsilon$.
    The completed relative classes sum to $A$ or $B$, so the total truncated exterior area is at most $2R$.
    Since $3(R-\epsilon)>2R$, there are at most two nonconstant exterior components.
    A closed exterior component would already have area at least $2R$, leaving no area for a punctured component connecting it to the bottom level.

    Due to the genus constraint, the component graph is a tree.
    Exactness and the noncontractibility of the flat Reeb orbits imply that every nonconstant component in the bottom or intermediate levels has at least two punctures.
    Every leaf is therefore exterior.
    There must be at least two leaves, so the tree is a path with exactly two exterior planes and all other components cylinders; in particular, there are no ordinary nodes or constant components.
    An intermediate cylinder with two negative ends would have $d\alpha_0$-energy equal to minus the sum of their Reeb actions, contradicting its nonnegativity \cite[Section~2.3]{DRGI}.
    Thus the path contains exactly one bottom cylinder, and each intermediate cylinder has one positive and one negative end.

    In class $A$, the cap intersections identify the relative classes of the two planes as $H_1^-$ and $H_1^+$.
    Their areas are $R$, and their boundary classes in \eqref{eq-boundary} determine the primitive asymptotic directions $\pm e_1$.
    By \cite[Lemma~4.2]{DRGI}, the bottom cylinder consequently has image $\{y=c,\ p_y=b\}$, and \cite[Lemma~4.4]{DRGI} identifies the intermediate cylinders as trivial.
    The class-$B$ case exchanges the coordinates.
    Each plane is simply covered because it meets one cap once.
    For limits of embedded ruling spheres, embeddedness then follows from the positivity-of-intersection argument in \cite[proof of Lemma~5.12]{DRGI}.
\end{proof}

\begin{remark}
    The above argument may replace part of the genericity argument in \cite[Corollary~5.4]{DRGI}, but does not generalize to the non-monotone case, to the best of the authors' knowledge.
    This allows us to keep the prescribed $\rho$-invariant almost complex structure $J_{\mathrm{out}}$ throughout the splitting argument without going into more complicated transversality issues, albeit achievable.
\end{remark}

We now use point constraints to obtain an exterior plane at every coordinate Reeb orbit.
Fix $z_0=(x_0,y_0,0,0)\in Z$ and, for each $\ell$, take the $J_\ell^{\mathrm{pre}}$-holomorphic ruling sphere in class $A$ through $z_0$.
SFT compactness with this point constraint gives a limiting building with a nonconstant bottom component through $z_0$, as in \cite[Proposition~5.7]{DRGI}.
By Lemma~\ref{lemma-splitting}, its bottom cylinder is $\{y=c,\ p_y=b\}$.
Since it contains $z_0$, we have $c=y_0$ and $b=0$.
For class $B$, the same argument gives $\{x=x_0,\ p_x=0\}$.
The two ends of each cylinder correspond to the two orientations of its coordinate geodesic.
As $y_0$ or $x_0$ varies, every signed coordinate Reeb orbit therefore occurs as the asymptote of an embedded exterior plane, paired with an opposite plane in the same limiting sphere.

\subsection{The four families of planes}\label{subsection-families}

For $c\in\R/\Z$, denote the unparametrized coordinate Reeb orbits in $S^*T^2$ by

\begin{equation}
    \gamma_1^\pm(c)={(x,c,\pm1,0)\in\R/\Z},\qquad
    \gamma_2^\pm(c)={(c,y,0,\pm1)\in\R/\Z}.
\end{equation}
Their orientations are given by $R_{\alpha_0}$.

\begin{figure}[htbp]
    \centering
    \includegraphics[width=9cm]{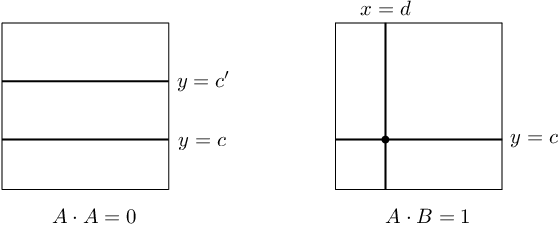}
    \caption{The standard inner cylinders used in Lemma~\ref{lemma-families}, shown at zero momentum.
        Cutting off the punctured ends and reconnecting matching boundary circles by annuli gives closed cycles with the displayed homological pairings.
        In each case the bottom intersections exhaust the pairing, so positivity forces all exterior intersections to vanish.}
    \label{fig-completion-intersections}
\end{figure}

\begin{lemma}\label{lemma-families}
    For every $i\in\{1,2\}$, each choice of sign, and every $c\in\R/\Z$, there is a unique exterior plane $P_i^\pm(c)$ asymptotic to $\gamma_i^\pm(c)$ that meets $D_i^\pm$ once and is disjoint from the other caps.
    For each $i$ and each choice of sign, these planes vary smoothly with $c$, and their union is an embedded three-dimensional submanifold of $Y\setminus Z$.
    Any two distinct planes among these four families are disjoint.
    Moreover,
    \begin{equation}\label{eq-plane-action}
        \rho P_1^\pm(c)=P_1^\pm(-c),\qquad
        \rho P_2^\pm(c)=P_2^\mp(c+\tfrac12).
    \end{equation}
\end{lemma}
\begin{proof}
    With $D_\infty=D_1^+\cup D_2^+$, our negative and positive cap labels correspond to the small and big planes of \cite{DRGI}, respectively.
    Lemma~\ref{lemma-splitting} and the preceding point-constrained construction supply the splitting configuration and embedded planes at every coordinate orbit.
    These planes have Maslov index $2$ and Fredholm index $1$, so the automatic transversality argument in \cite[proof of Proposition~5.11(i)]{DRGI} applies to our fixed $J_{\mathrm{out}}$.
    The arguments of \cite[Proposition~5.11 and Lemmas~5.12--5.13]{DRGI} therefore give uniqueness, smooth parametrization by $c$, and an embedded three-dimensional union for each family.

    Any two distinct planes in the same coordinate direction are disjoint by \cite[Lemmas~5.5 and~5.8]{DRGI}.
    For different coordinate directions, complete the planes to split spheres in classes $A$ and $B$.
    Their bottom cylinders $\{y=c,p_y=0\}$ and $\{x=d,p_x=0\}$ intersect once.
    The completion comparison in \cite[proof of Lemma~5.8]{DRGI} then accounts for $A\cdot B=1$, so positivity excludes any exterior intersection.
    See Figure~\ref{fig-completion-intersections}.

    Finally, $\rho$ preserves $J_{\mathrm{out}}$ and agrees with $\sigma$ near $D\cup Z$.
    It permutes the cap and orbit data as prescribed in \eqref{eq-plane-action}, so uniqueness gives that formula.
\end{proof}

Each of the four families $P_i^\pm$ therefore fills a hypersurface of $Y\setminus Z$.
In the next subsection, we recover the two sphere rulings by joining the opposite families.

\subsection{Smoothing and the two sphere rulings}\label{subsection-smoothing}

\begin{proof}[Proof of Proposition~\ref{prop-marking}]
    Apply the end-straightening construction of \cite[Proposition~5.16]{DRGI} to the four families in Lemma~\ref{lemma-families}.
    This modifies the planes near $Z$, preserving their symplecticity and the smooth embedded families, so that their ends become standard half-cylinders.
    The construction can be performed simultaneously and $\rho$-equivariantly: the four families of ends lie in disjoint cones around the momentum directions $\pm e_1,\pm e_2$, and we use the same radial cutoff in all four cones.
    Indeed, $\rho=\sigma$ near $Z$ acts by translations and sign changes, permutes these cones, and preserves the radial coordinate, so the cutoff construction commutes with $\rho$.
    See Figure~\ref{fig-simultaneous-smoothing}.

    For each $c$, the straightened planes $P_1^\pm(c)$ join smoothly across $\{y=c,p=0\}$ to form a symplectic sphere in class $A$, with local model $\{y=c,p_y=0\}$ near $Z$.
    Similarly, the planes $P_2^\pm(d)$ join to a sphere in class $B$, with local model $\{x=d,p_x=0\}$.
    These give two smooth circle families of spheres.
    Within each family the spheres are disjoint, while spheres from different families meet only at their coordinate-circle crossing on $Z$.
    Near $Z$, their tangent planes are the two standard, symplectically orthogonal coordinate planes.
    Both families are preserved by $\rho$, and the construction leaves the planes unchanged near $D$.

    For a single sphere family, the existence of a global tame almost complex structure making all its leaves holomorphic is noted in \cite[Proposition~5.16(iii) and Remark~5.17]{DRGI}.
    Here we need one compatible, $\rho$-invariant almost complex structure for both families.
    Each family foliates an embedded three-dimensional submanifold, along which the sphere tangent planes form a smooth symplectic subbundle $V$ of $TY$.
    Choose a Riemannian metric $g$ satisfying $V\perp_g V^\omega$ along every leaf.
    Near $Z$, the standard product metric satisfies this condition for both families.
    Near $D$, we let the metric be associated to $J_{\mathrm{out}}=J_{\mathrm{std}}$.
    Away from a neighborhood of $Z$, the two family unions are disjoint, so these metric prescriptions extend in separate tubular neighborhoods and then to $Y$.
    A partition of unity and averaging over $\rho$ preserve the orthogonality conditions and the prescribed local metrics.
    The compatible almost complex structure $J_1$ associated to this metric and $\omega$ therefore preserves every sphere tangent plane, is $\rho$-invariant, and is standard near $D\cup Z$.
    Gromov's theorem identifies the constructed spheres as leaves of its two sphere rulings, proving the proposition.

    \begin{figure}[htbp]
        \centering
        \includegraphics[width=13cm]{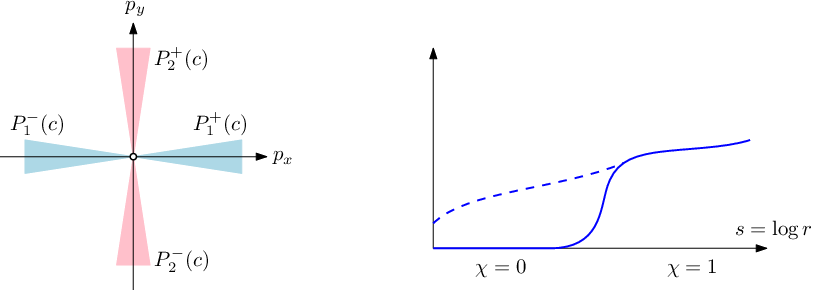}
        \caption{Pairing the ends after simultaneous smoothing in Proposition~\ref{prop-marking}.
            Left: near $Z$, the ends of $P_i^\pm(c)$ lie in four disjoint cones of nonzero momenta.
            Base coordinates are suppressed.
            Right: one scalar end-graph component at the fixed $(x,c)$, before and after multiplication by $\chi(s)$.
            The modified graphs vanish for $s\ll0$ and agree with the original ones near the outer boundary.
            The same scalar cutoff in all four cones commutes with $\rho$.
            The standard opposite ends then close across their common geodesic on $Z$.}
        \label{fig-simultaneous-smoothing}
    \end{figure}
\end{proof}

\section{Relative reconstruction}\label{section-reconstruction}

We prove Proposition~\ref{prop-equivariant} using Gromov's ruling argument, as in \cite[Section~5.5]{DRGI}, with the reconstruction described in \cite[pp.~11--13]{Evans}.
We arrange that the maps conjugate $\rho$ to $\sigma$, preserve each labeled cap, and preserve $Z$ at the endpoint.

\subsection{Leaf intersection maps}\label{subsection-leaf-maps}

Set $J_0=\psi_*J_{\mathrm{std}}$.
Choose a smooth path of $\rho$-invariant $\omega$-compatible almost complex structures $J_t$ from $J_0$ to $J_1$, equal to $J_{\mathrm{std}}$ on one fixed neighborhood of $D$.
Such a path is obtained by interpolating the associated invariant metrics, as in \cite[Proposition~12.8]{CdS}.

Let $A_t(p)$ and $B_t(p)$ denote the two ruling leaves through $p\in Y$.
We use the pole labels and factor involutions introduced in Subsection~\ref{section-square}.
Parametrize $D_1^-$ by $v(b)=(a_-,b)$.
Choose a $\sigma_2$-fixed point $b_0$ on the second equator and set $P=(a_-,b_0)$.
Since $\rho=\sigma$ near $D$, the point $P$ is fixed by $\rho$; it lies on the equator of $D_1^-$, away from the corners $(a_-,b_\pm)$.
Uniqueness of the ruling leaf through $P$ makes $A_t(P)$ $\rho$-invariant.
Its unique intersection with $D_1^+$ is therefore fixed by $\rho$.
This intersection varies continuously in the two-point fixed set of $\rho|_{D_1^+}$ and starts at $(a_+,b_0)$, so it remains there.
Thus $A_t(P)$ meets $D_1^\pm$ at $(a_\pm,b_0)$ for every $t$.

Choose smooth orientation-preserving parametrizations $u_t:S^2\to A_t(P)$ satisfying
\[
    u_0(a)=\psi(a,b_0),\qquad u_t\sigma_1=\rho u_t,
\]
with first coordinate $a$ near $a_\pm$.
To obtain them, use the equivariant holomorphic graph parametrizations near the caps, extend the normal velocity of $A_t(P)$ with these prescribed graph velocities, and average the extension over $\rho$.
Its flow transports $u_0$ to the required $u_t$.

\begin{figure}[htbp]
    \centering
    \includegraphics[width=7cm]{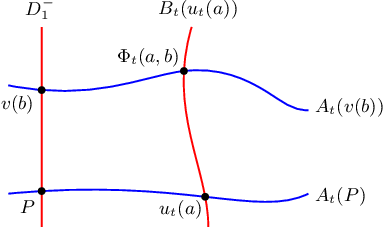}
    \caption{Reconstruction from the two sphere rulings.
        Lines denote sphere leaves.
        The point $\Phi_t(a,b)$ is the intersection of the $B$-leaf through $u_t(a)$ and the $A$-leaf through $v(b)$.}
    \label{fig-leaf-reconstruction}
\end{figure}

Define
\begin{equation}\label{eq-phi}
    \Phi_t:S^2\times S^2\to Y,\qquad
    \Phi_t(a,b)=B_t(u_t(a))\cap A_t(v(b)).
\end{equation}
The leaf-intersection construction of \cite[p.~11]{Evans} gives a smooth family of diffeomorphisms.
The chosen parametrizations, together with $v\sigma_2=\rho v$, give
\begin{equation}\label{eq-phi-identities}
    \Phi_0=\psi,\qquad \Phi_t\sigma=\rho\Phi_t.
\end{equation}
The pole labels show that each $\Phi_t$ preserves every labeled cap.
Moreover, product spheres sufficiently close to the caps are $J_t$-holomorphic and hence are ruling leaves, as in \cite[proof of Theorem~D, Section~5.4]{DRGI}.
The chosen parametrizations label these leaves by their product coordinates, so $\Phi_t$ is the identity on one common neighborhood of the four corners.
Figure~\ref{fig-leaf-reconstruction} illustrates \eqref{eq-phi}.

\begin{figure}[htbp]
    \centering
    \includegraphics[width=9cm]{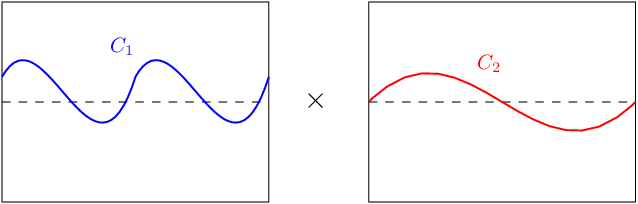}
    \caption{Before the appended straightening, the inverse leaf map sends $Z$ to a product of invariant essential circles, drawn schematically.
        The appended equivariant annular isotopies take the circles to the equators and hence their product to $Z$.}
    \label{fig-circle}
\end{figure}

Set $f_t=\Phi_t^{-1}$.
By Proposition~\ref{prop-marking},
\[
    f_1(Z)=C_1\times C_2,
\]
where $C_i$ is a $\sigma_i$-invariant embedded circle in the corresponding annular factor.
Both circles are essential because $Z\hookrightarrow Y\setminus D$ induces an isomorphism on fundamental groups and $f_1$ preserves $D$.
Apply Lemma~\ref{lemma-annulus} in each factor and append the resulting product isotopy; see Figure~\ref{fig-circle}.
After stationary reparametrizations at the join, the resulting smooth path, still denoted $f_t$, satisfies
\begin{equation}\label{eq-reconstruction-path}
    f_0=\psi^{-1},\qquad f_t\rho=\sigma f_t,\qquad
    f_t(D_i^\pm)=D_i^\pm,\qquad f_1(Z)=Z.
\end{equation}
It remains the identity near the corners.
Throughout the original interval, $f_t$ takes the two $J_t$-ruling planes to the product-factor planes and preserves their orientations.
This remains true on the appended interval with the almost complex structure fixed at $J_1$.

\subsection{Relative Moser deformation and cutoff}\label{subsection-moser}\label{subsection-cutoff}

Consider
\begin{equation}\label{eq-moser-forms}
    \omega_{s,t}=(1-s)\omega+s f_t^*\omega,\qquad 0\leq s,t\leq1.
\end{equation}
These forms are symplectic by the convex-form argument of \cite[pp.~12--13]{Evans}: the two ruling planes are symplectically orthogonal with respect to $f_t^*\omega$, which is positive on each and hence tames the same $J_t$ as $\omega$.
On the appended interval in the previous subsection, use $J_1$.
The forms are $\rho$-invariant by \eqref{eq-reconstruction-path} and cohomologous because $f_t$ starts at the symplectomorphism $\psi^{-1}$.

We next apply a relative equivariant Moser argument to make the path $f_t$ symplectic, while maintaining its initial value and zero-section endpoint condition.
A subsequent Hamiltonian cutoff makes the corrected path equal to the identity on a fixed neighborhood of $D$.

To keep the zero section at the endpoint, we need
\begin{equation}\label{eq-relative}
    [f_1^*\omega-\omega]=0\quad\text{in }H^2(Y,Z;\R).
\end{equation}
The map $f_1$ preserves $Z$ and every labeled cap.
It also preserves each cap's normal orientation, since these signs are constant along the cap-preserving path starting at $\psi^{-1}$.
It therefore preserves the cap-intersection coordinates of $H_2(Y,Z;\Z)$.
The hemisphere basis in Subsection~\ref{section-square} shows that $f_1$ acts trivially on this group, proving \eqref{eq-relative}.

The homotopy formula for $f_t$ gives smooth primitives $\eta_t$ of $f_t^*\omega-\omega$ with $\eta_0=0$.
Average them over $\rho$, and choose a $\rho$-invariant primitive $\beta^{\mathrm{rel}}$ at $t=1$ with $\beta^{\mathrm{rel}}|_{TZ}=0$, using \eqref{eq-relative}.
Since $\beta^{\mathrm{rel}}-\eta_1$ is closed, setting
\[
    \beta_t=\eta_t+t(\beta^{\mathrm{rel}}-\eta_1)
\]
gives invariant primitives satisfying
\begin{equation}\label{eq-primitive-endpoints}
    \beta_0=0,\qquad \beta_1|_{TZ}=0.
\end{equation}

Near the corners, $f_t=\id$, so $\beta_t$ is closed.
Write $\beta_t=db_t$ on disjoint corner balls, with $b_t$ normalized to vanish at the corners.
Subtracting $d(\chi b_t)$ with invariant cutoffs makes $\beta_t$ vanish on smaller common corner neighborhoods.
This preserves \eqref{eq-primitive-endpoints}, since the balls are disjoint from $Z$ and the normalized potentials vanish at $t=0$.

Each cap $C$ is symplectic for $\omega_{s,t}$.
Choose functions $a_{s,t}^C$ vanishing on $C$ whose normal derivatives satisfy
\begin{equation}\label{eq-cap-jet}
    da_{s,t}^C|_{(TC)^{\omega_{s,t}}}
    =-\beta_t|_{(TC)^{\omega_{s,t}}}.
\end{equation}
A tubular extension and cutoff realize these derivatives.
The extensions can be chosen smoothly and linearly in $\beta_t$, with disjoint supports away from the corner neighborhoods and $Z$.
Let $a_{s,t}$ be the sum of these functions averaged over $\rho$, and define $\widehat\beta_{s,t}=\beta_t+da_{s,t}$.
Averaging preserves the prescribed derivatives.
The corrected primitives vanish at $t=0$, annihilate $TZ$ at $t=1$, and annihilate $(TC)^{\omega_{s,t}}$ along every cap.

The Moser vector field defined by $\iota_{V_{s,t}}\omega_{s,t}=-\widehat\beta_{s,t}$ is therefore $\rho$-invariant, tangent to every cap, zero at $t=0$, and tangent to $Z$ at $t=1$.
Its flow from $s=0$ to $s=1$, denoted by $M_t$, satisfies $M_t^*f_t^*\omega=\omega$.
Thus
\begin{equation}\label{eq-symplectic-path}
    S_t=f_tM_t
\end{equation}
is a smooth symplectic path with
\[
    S_t\rho=\sigma S_t,\qquad S_0=\psi^{-1},\qquad
    S_t(D_i^\pm)=D_i^\pm,\qquad S_1(Z)=Z.
\]

We now make $S_t$ equal to the identity near $D$ by a Hamiltonian cutoff, following the construction at the end of \cite[Section~7.2]{DRGI}.
Since $H^1(Y;\R)=0$, the vector field $X_t=\dot S_tS_t^{-1}$ has a smooth Hamiltonian $H_t$ with $\iota_{X_t}\omega=-dH_t$.
The conjugacy relation makes $X_t$ $\sigma$-invariant, so averaging allows us to take $H_t$ invariant as well.

Choose an invariant cutoff $\chi$ supported near $D$, away from $Z$, and equal to one on a neighborhood $O$ of $D$.
Let $K_t$ be the Hamiltonian flow of $\chi H_t$ from the identity.
It commutes with $\sigma$ and fixes a neighborhood of $Z$.
Because $S_t(D)=D$ and the parameter interval is compact, there is an invariant neighborhood $N$ of $D$ on which $\psi=\id$ and such that $S_t(N)\subset O$ for all $t$.
On these trajectories the two Hamiltonian vector fields agree, and their initial maps agree on $N$.
Uniqueness of the flow gives
\begin{equation}\label{eq-equal-cap-germs}
    K_t|_N=S_t|_N\qquad\text{for all }t.
\end{equation}

Define
\begin{equation}\label{eq-equivariant-path}
    h_t=K_t^{-1}S_t\psi.
\end{equation}
This is a symplectic isotopy from the identity, commuting with $\sigma$ and equal to the identity on $N$.
Its supports lie in the common compact set $Y\setminus N$.
Restricting to $Y\setminus D$ and extending by the identity gives the required isotopy on $T^*T^2$.
Finally,
\[
    h_1\bigl(\psi^{-1}(Z)\bigr)=K_1^{-1}S_1(Z)=Z,
\]
because $K_1$ fixes $Z$.
This proves Proposition~\ref{prop-equivariant}.

\section{Conclusion of the proof}\label{section-descent}

\begin{proof}[Proof of Theorem~\ref{thm-main}]
    Let $E$ and $\psi$ be as in Lemma~\ref{lemma-normalization}.
    Proposition~\ref{prop-equivariant} gives a smooth symplectic isotopy $h_t$ from the identity with $h_t\sigma=\sigma h_t$ and $h_1(E)=Z$.
    For the covering projection $\pi:T^*T^2\to T^*K$, equivariance gives a smooth isotopy $\overline{h}_t$ defined by
    \begin{equation*}
        \overline{h}_t\pi=\pi h_t,
    \end{equation*}
    which is symplectic since $\pi$ is locally symplectic.
    Its initial map is the identity and $\overline{h}_1(L)$ is the zero section.
    If $C\subset T^*T^2$ is a common compact support of $h_t$, then $\pi(C)$ is a common compact support for $\overline{h}_t$.

    Since $T^*K$ is an oriented four-manifold and retracts onto $K$, Poincar\'e duality with compact supports gives
    \begin{equation*}
        H_c^1(T^*K;\R)\cong H_3(T^*K;\R)\cong H_3(K;\R)=0.
    \end{equation*}
    The descended path is therefore generated by a smooth Hamiltonian with one common compact support, by the Hamiltonian isotopy criterion in \cite[Proposition~9.3.1 and p.~389]{MS2017}.
    This proves the theorem.
\end{proof}

\appendix
\section{Equivariant annular straightening}\label{appendix-annulus}

In this appendix we prove the equivariant straightening result used to complete the smooth reconstruction in Subsection~\ref{subsection-leaf-maps}.
The inverse leaf-intersection map at the endpoint sends $Z$ to a product $C_1\times C_2$ of essential circles invariant under $\sigma_1$ and $\sigma_2$, respectively.
The main result is Lemma~\ref{lemma-annulus}, which moves each circle to the equator through an isotopy commuting with the corresponding involution and fixed near the annular ends.
Appending the product of these isotopies gives $f_1(Z)=Z$ while preserving the conjugacy relation and every labeled cap, providing the endpoint condition needed for the relative Moser deformation in the proof of Proposition~\ref{prop-equivariant}.

\begin{lemma}\label{lemma-annulus}
    On $\mathcal A=S^1\times\R$, let $S(q,p)=(q+\tfrac12,p)$ and $r(q,p)=(-q,-p)$.
    If $C$ is a smooth invariant embedded essential circle for either involution $\tau=S$ or $\tau=r$, then there is a smooth isotopy $f_t$ with one common compact support such that
    \begin{equation*}
        f_0=\id,\qquad f_t\tau=\tau f_t,\qquad f_1(C)=S^1\times\{0\}.
    \end{equation*}
\end{lemma}

\begin{proof}
    For $S$, the quotient is an annulus and the image of $C$ is an embedded essential circle.
    Curve isotopy and smooth isotopy extension \cite[Propositions~1.10--1.11]{FM} give an ambient isotopy supported in a compact subannulus carrying it to the core.
    Lifting the vector field gives the required isotopy upstairs.

    For $r$, choose a large symmetric closed annulus containing $C$.
    Its quotient is a disk with two simple branch points, locally described by $w=z^2$.
    See Figure~\ref{fig-annular-branch-lift}.
    Since $r$ acts by $-1$ on the first homology of the annulus, its restriction to $C$ reverses orientation, hence $C$ contains the two fixed points and projects to an embedded arc joining the branch points.
    We first straighten $C$ near those points by equivariant local isotopies.
    After rotating its tangent line, write it as $y=h(x)$ with $h$ odd.
    The vector field $-\chi(x,y)h(x)\partial_y$, with an even cutoff equal to one on the local trace, straightens this graph.
    The two constructions have disjoint supports.

    The coordinate-change principle for surfaces \cite[Section~1.3]{FM} gives a disk diffeomorphism, fixed near the boundary, taking the resulting arc to the image of the core.
    Since the arcs have straight endpoint germs, we specify complex linear germs in the centered branch charts.
    The mapping class group of the disk with two unordered marked points is $B_2\cong\Z$ \cite[Sections~9.1--9.2]{FM}.
    A generating half twist preserves the target arc and can also have complex linear branch germs.
    Postcomposing with a suitable power gives a straightening map $G$ isotopic to the identity relative to the boundary and the marked points.
    Choose a smooth isotopy $G_t$ from the identity to $G$, fixed on one boundary collar; smooth isotopies may be used by \cite[discussion after Theorem~1.13]{FM}.
    Figure~\ref{fig-annular-branch-lift} illustrates the branched quotient and the straightening of the projected arc.

    \begin{figure}[htbp]
        \centering
        \includegraphics[width=7cm]{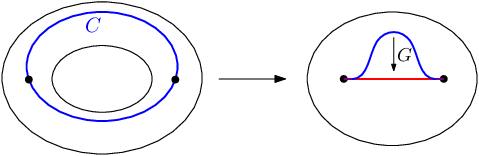}
        \caption{The branched case of Lemma~\ref{lemma-annulus}.
            After the initial local straightening, the invariant circle descends to an arc joining the two branch points.
            The disk isotopy after the half-twist correction takes this arc to the image of the core.}
        \label{fig-annular-branch-lift}
    \end{figure}

    Before lifting, we normalize the branch germs relative to the endpoints of the isotopy.
    In a centered branch chart, we write
    \begin{equation*}
        G_t(w)=A_tw+R_t(w),\qquad
        A_t\in\operatorname{GL}^+(2,\R),\qquad R_t(w)=O(|w|^2),
    \end{equation*}
    uniformly in $t$.
    For sufficiently small fixed $\epsilon$, replace this expression by $G_t(w)-\chi(w/\epsilon)R_t(w)$, which is exactly $A_tw$ on a smaller disk.
    The change is $O(\epsilon)$ in $C^1$, so the maps remain diffeomorphisms.
    The endpoints are unchanged since their germs are already linear.

    We write $A_t=U_tP_t$ in polar form and define
    \begin{equation*}
        B_t=P_t^{-1}\sqrt{\det P_t}\,I.
    \end{equation*}
    Precompose with the time-$1$ flow of the compactly supported vector field $\chi_0(w)(\log B_t)w$.
    On a sufficiently small common disk this flow has germ $B_t$, making the new germ $U_t\sqrt{\det P_t}\,I$ which is complex-linear.
    The correction is the identity at both endpoints since $B_0=B_1=I$.
    We then perform these constructions in the two disjoint branch charts.

    The resulting germs are $w\mapsto a_tw$ with $a_t\ne0$ and $a_0=1$.
    At each branch point, choose the smooth square root along the parameter interval that starts at $1$.
    In the branch chart $w=z^2$, the local lift $\widetilde G_t(z)=\sqrt{a_t}\,z$ is smooth in $(t,z)$ and satisfies
    \begin{equation*}
        (\widetilde G_t(z))^2=a_tz^2=G_t(z^2).
    \end{equation*}
    These local lifts agree with the ordinary isotopy lift away from the branch points by uniqueness and their common initial value.
    The lift starting at the identity commutes with $r$ and is fixed near both annular ends.
    We concatenate it with the initial local straightening to get the required isotopy.
    All constructions take place in one compact annulus, and stationary reparametrizations make the concatenation smooth.
\end{proof}

\bibliographystyle{alpha}
\bibliography{klein}

\end{document}